\documentclass{article}

\usepackage[utf8]{inputenc}
\newcounter{defn}
\newcounter{rmk}

\usepackage{amsmath}
\usepackage{float}

\usepackage{comment}
\usepackage[linktocpage=true]{hyperref}
\usepackage[margin=1.2 in, top=1 in, bottom= 1.2 in]{geometry}

\usepackage{indentfirst} 
\usepackage[english]{babel}
\usepackage{amsfonts}
\usepackage{mathrsfs}
\usepackage[mathscr]{euscript}
\usepackage{bbm}
\usepackage{latexsym}
\usepackage{mathdots}
\usepackage{amssymb}
\usepackage{mathtools}
\usepackage{graphicx}
\usepackage{tikz}
\usepackage{placeins}
\usepackage{esint}

\usepackage{longtable}
\usepackage{booktabs}
\usepackage{array}

\newcolumntype{L}{>{\raggedright\arraybackslash}p{0.22\textwidth}}
\newcolumntype{D}{>{\raggedright\arraybackslash}p{0.58\textwidth}}
\newcolumntype{R}{>{\raggedright\arraybackslash}p{0.12\textwidth}}

\usepackage{enumitem}
\setlist{nolistsep}
\usepackage{amsthm}
\usepackage{relsize}
\usepackage{nicefrac}
\usepackage[capitalize]{cleveref}

\newtheoremstyle{plain}{3mm}{3mm}{\slshape}{}{\bfseries}{.}{.5em}{}
\newtheoremstyle{definition}{2mm}{2mm}{}{}{\bfseries}{.}{.5em}{}
\theoremstyle{plain}
	
\newtheorem{theorem}{Theorem}

\newtheorem{lemma}[theorem]{Lemma}

\theoremstyle{definition}
\newtheorem{definition}[defn]{Definition}
\newtheorem{remark}[rmk]{Remark}

\theoremstyle{plain}
\newtheorem*{namedthm}{\namedthmname}
\newcounter{namedthm}
	

\newcommand{\R}{\mathbb{R}}

\newcommand{\Gee}{\mathbb{G}}
\newcommand{\diam}{\text{diam}}

\newcommand{\F}{\mathbb{F}}

\newcommand{\eps}{\epsilon}

\newcommand{\B}{\mathcal{B}}

\newcommand{\Tjj}{T_{Q^{1}_{j}}}

\newcommand{\ree}{\mathbb{R}^{n+1}}
\newcommand{\re}{\mathbb{R}}
\title{BMO solvability with singular drifts on ample sawtooth domains implies $L^p$ solvability.}
\author{Aritro Pathak}

\date{}

\begin{document}
\maketitle
\begin{abstract}
  For a linear elliptic operator with a singular drift that satisfies a finite Carleson measure condition, we prove that there exist `ample' sawtooth domains of the unit ball $B(0,1)\subset \R^{n+1}$ so that a BMO solvability assumption in these sawtooth subdomains implies that the elliptic measure satisfies the weak $A_\infty$ condition with respect to the surface measure on this `ample' sawtooth domain. This is a quantifiable absolute continuity condition, which is equivalent to saying the $L^p$ Dirichlet problem is solvable for some $1<p<\infty$. Such singular drifts have been considered in the literature in the context of perturbative $L^p$ Dirichlet solvability problems, by Hofmann-Lewis and Kenig-Pipher. By an ample sawtooth domain, we mean a sawtooth domain whose boundary coincides with the boundary of the unit ball, except for an arbitrarily small fraction. The methods can be naturally extended to show the result for more general bounded Lipschitz domains.
  \end{abstract}

\maketitle

\section{Introduction}

There is an extensive literature on the connection between $L^p$ solvability of the Dirichlet problem, and scale invariant quantifiable absolutely continuity property of the elliptic measures corresponding to divergence form elliptic equations. Further work of \cite{DKP11, Z18,HL18} explored the connection between a natural BMO estimate on the solutions of the Dirichlet problem with continuous data, and the quantifiable absolute continuity condition of the elliptic measure and thus also the question of $L^p$ solvability.

Our goal is to understand how the presence of a singular lower order drift term, satisfying a Carleson measure condition, affects the relationship between BMO-solvability and the weak-$A_\infty$ 
property of the elliptic measure. In this manuscript, we assume this natural $BMO$ solvability estimate on the solutions, for the operator of \cref{operator} , which has a divergence form principal term along with a lower order singular drift term. Such drift terms satisfying a natural Carleson measure condition, have been considered previously for perturbative questions for Dirichlet solvability \cite{HL01, KP01}.

In particular, we show that the mechanism developed in \cite{HL18} continues to hold in this perturbed setting for a bounded domain such as the unit ball, by introducing significant extensions using the Markov property for the elliptic measure.

Consider the unit ball $B(0,1)\subset \R^{n+1}$, $n \geq 2$,  and the operator,
\begin{align}\label{operator}
    L:=-\text{div}(A(x)\nabla (\cdot)) + \B\cdot\nabla(\cdot).
\end{align}

Here $A(x)$ is an elliptic matrix, so that, we have parameters $0<\lambda<\Lambda<\infty $, so that,
\begin{align}
    |A_{ij}(x)|\leq \Lambda, \text{for each}\  0\leq i,j\leq n+1, A_{ij}(x)\eta_i\eta_j\geq \lambda|\eta|^2,
\end{align}
uniformly over the domain $B(0,1)$, where $\eta$ is a row vector.

Here, $\partial\Omega=S(0,1)$ is the boundary of the domain. The drift $\B$ satisfies a Carleson measure condition, which for the unit ball setting, is the requirement that, 
\begin{align}
    \label{DKP}
    \sup\limits_{x\in \partial\Omega,0<r<1/2}\frac{1}{\sigma(\Delta(x,r))}\int\limits_{\Omega\cap B(x,r)}\sup\limits_{y\in B(t,\delta(t)/2)}\Big(|\B|^{2}(y)\delta(y)\Big) dV(t)=\mathscr{M}<\infty.
\end{align}

Here $\delta_{\Omega}(x):=\text{dist}(x,\partial\Omega)$ denotes the distance to the boundary of the domain $\Omega$.

Note that the supremum is over all balls centered on the boundary, and of radius less than or equal to $1/2$, and thus the union of these balls is disjoint from $B(0,1/2)$.

By hypothesis we also have some $M<\infty$, such that we have the ambient pointwise upper bound condition,
\begin{align}\label{large}
    |\B(x)|\leq \frac{M}{\delta_{\Omega}(x)},
\end{align}
for any $x\in \R^{n+1} $.

We note that the non-divergence form elliptic operator can be written as a divergence form operator along with a drift term, as in \cref{operator}, and the nondivergence form problem is one of the primary motivations for studying the Dirichlet solvability question with the drift term\cite{KP01,HL01}.

\subsection{Preliminaries and statement of result.}

We recall standard notions such as the dyadic decomposition of Ahlfors regular sets, the Whitney decomposition of the unit ball, and the definitions of BMO-solvability and weak-$A_\infty$.

\begin{lemma}\label{lemmaCh}({Existence and properties of the ``dyadic grid''})
\cite{DS1,DS2}, \cite{Ch90}.
Suppose that $E\subset \ree$ is a closed $n$-dimensional ADR set.  Then there exist
constants $ a>0,\, \gamma>0$ and $C_*<\infty$, depending only on dimension and the
ADR constant, such that for each $k \in \mathbb{Z},$
there is a collection of Borel sets (``cubes'')
$$
\mathbb{D}_k:=\{Q_{j}^k\subset E: j\in \mathfrak{I}_k\},$$ where
$\mathfrak{I}_k$ denotes some (possibly finite) index set depending on $k$, satisfying

\begin{list}{$(\theenumi)$}{\usecounter{enumi}\leftmargin=.8cm
\labelwidth=.8cm\itemsep=0.2cm\topsep=.1cm
\renewcommand{\theenumi}{\roman{enumi}}}

\item $E=\cup_{j}Q_{j}^k\,\,$ for each
$k\in{\mathbb Z}$.

\item If $m\geq k$ then either $Q_{i}^{m}\subset Q_{j}^{k}$ or
$Q_{i}^{m}\cap Q_{j}^{k}=\emptyset$.

\item For each $(j,k)$ and each $m<k$, there is a unique
$i$ such that $Q_{j}^k\subset Q_{i}^m$.

\item $\diam\big(Q_{j}^k\big)\leq C_* 2^{-k}$.

\item Each $Q_{j}^k$ contains some ``surface ball'' $\Delta \big(x^k_{j},a2^{-k}\big):=
B\big(x^k_{j},a2^{-k}\big)\cap E$.

\item $H^n\big(\big\{x\in Q^k_j:{\rm dist}(x,E\setminus Q^k_j)\leq \varrho \,2^{-k}\big\}\big)\leq
C_*\,\varrho^\gamma\,H^n\big(Q^k_j\big),$ for all $k,j$ and for all $\varrho\in (0,a)$.
\end{list}
\end{lemma}

This result above is stated in much larger generality for any domain that is the complement of an Ahlfors-David regular set. For the unit ball under consideration, we write $\mathbb{D}(S(0,1))=\sqcup_{k\geq 2} \mathbb{D}_k $ and any cube belonging to $\mathbb{D}_k$ has length comparable to $2^{-k}$ using the symmetry of the unit ball. We denote the surface measure of the cube $Q^{k}_{j}\subset S(0,1)$ by $H^{n}(Q^{k}_{j}):=\sigma(Q_{j}^{k})$.

Corresponding to the dyadic decomposition of Lemma 1, we have the Whitney decomposition of the unit ball $\mathcal{W}(B(0,1))$. 
For our specific case of the unit ball, it is enough to consider that for some fixed $k$, the cubes belonging to $\mathbb{D}_k$ are disjoint.  Note that in this process, we leave out the the ball $B(0,1/2)$ from the Whitney decomposition since we have considered only the union of the sets $\mathbb{D}_k$ for $k\geq 2$.

By the symmetry of the unit ball, we see that for any $U_{Q}\in \mathcal{W}(B(0,1))$, with $Q\in \mathbb{D}_{k+1}$, we have, $U_{Q}=Q\times \{r:1- 2^{-k}\leq r<1-2^{-k-1}\}, T_Q= Q\times \{r:1-2^{-k}<r<1\} $. Further, note that the radial width of $U_Q$ with $Q\in \mathbb{D}_k$, is $l(Q)=2^{-(k+1)}$, which we also call the length of $Q\in \mathbb{D}_{k+1}$ as well as the length of $U_Q$.


Note that the Carleson measure condition of \cref{DKP} is equivalent to the condition in terms of the Carleson boxes, that,
\begin{align}
    \label{DKP2}
    \sup\limits_{Q\in \mathbb{D}(S(0,1))}\frac{1}{\sigma(Q)}\int\limits_{T_Q}\sup\limits_{y\in B(t,\delta(t)/2)}\Big(|\B|^{2}(y)\delta(y)\Big) dV(t)<\infty.
\end{align}

Fix an $\eps$ small enough, to be specified later. We define the integer $l_0$ so that,



\begin{equation}
    \frac{M}{2^{-k}}=\overline{c}\frac{\eps}{\tau_{k}}, \text{with} \ \tau_{k}l_0= 2^{-k},
\end{equation}

and thus, 
\begin{equation}\label{eq5}
    l_{0}= \frac{M}{\overline{c}\eps}, \tau_k =\frac{\overline{c}\eps}{M }2^{-k}.
\end{equation}

The parameter $\overline{c}$ is chosen so that given $M$ and $\eps$, we have that $l_0=\lfloor M/\overline{c}\eps \rfloor$ is an integer. 

For each Carleson box $T_Q$ corresponding to the the Whitney cube $Q$, we write 
\begin{equation}
    T_{Q}=S_{Q}\sqcup U_Q,
\end{equation}

where $S_Q=T_{Q}\setminus U_Q$ is the `bottom half' of the Carleson box.




We make a refinement of $\mathcal{W}(B(0,1))$ so that each Whitney cube of this original decomposition is written as a disjoint union of
$(\lfloor2\sqrt{n+1}\rfloor)^{n+1}$ many congruent copies of the original cube, by taking these smaller cubes to have length reduced by a factor of $1/\lfloor2\sqrt{n+1}\rfloor$, so that the diagonal of these smaller cubes is $1/2$ the length of the original cube . We call this modified decomposition $\mathcal{W}'(B(0,1))$. 

It is clear that each Whitney cube in $\mathcal{W}'(B(0,1))$ is contained in some unique Whitney cube $U_Q \in \mathcal{W}(B(0,1))$. When we are restricted to some Carleson box $T_Q$, we denote by $\mathcal{W}(T_Q), \mathcal{W}'(T_Q)$ respectively the set of Whitney cubes of $\mathcal{W}(B(0,1)), \mathcal{W}'(B(0,1))$ restricted to $T_{Q}$.

For any Whitney cube $P\in \mathcal{W}' (B(0,1))$,  $P\subset B(0,1)$, we naturally denote by $l(P)$ the length of the cube. For any dyadic cube $Q\in \mathbb{D}$, $Q\subset S(0,1)$, we also use $l(Q)\approx l(U_Q)$ to denote the diameter of the dyadic cube up to a constant.

\begin{definition}[Average smallness assumption(ASA)]\label{ASA}
    Given $\eps>0$, a cube $P $ of the modified Whitney decomposition $\mathcal{W}'(B(0,1))$ satisfies the average smallness assumption, if we have,
    \begin{align}\label{avg}
        \int_{P} \ \sup\limits_{B(t,\delta(t)/2)}\Big(|\B|^{2}(y)\delta(y)\Big)dV(t) \leq \eps l(P)^{n}.
    \end{align}
    If the reverse inequality holds, the average smallness assumption fails.
\end{definition}

We note that the supremum condition within the integrand in \cref{DKP} means that an average smallness assumption (ASA) on a Whitney cube belonging in $\mathcal{W}'(B(0,1))$ implies pointwise smallness on that Whitney cube. For the cubes that satisfy the average smallness assumption(ASA), we will have with some $c$ constant uniform over the half space, that,
\begin{align}
    c\cdot \sup\limits_{P}\Big(|\B|^{2}(y)\delta(y)\Big) l(P)^{n+1}  \leq \int_{P} \ \sup\limits_{B(t,\delta(t)/2)}\Big(|\B|^{2}(y)\delta(y)\Big)dV(t) \leq \eps l(P)^{n}.
\end{align}

The left side inequality above is a consequence of the fact that for any point $t\in P$, the ball $B(t,\delta(t)/2)\supset P$ which in turn relies on the way the modified set $\mathcal{W}'(B(0,1))$ has been constructed, and the fact that the Whitney cubes are constructed as $U_Q=Q\times \{ r:1- 2^-k <r <1-2^{-k-1} \}$ with $Q\in \mathbb{D}_{k+1} $.

From here, the pointwise smallness of the drift over the ASA cube follows for arbitrarily $\eps$.
\begin{align}\label{smallness}
    |\B(x)|\leq \frac{c_1 \sqrt{\eps}}{(l(P))}, \ \text{when} \ x\in P.
\end{align}




\begin{definition}
    A cone of aperture $a$ for $Q\in \partial \Omega_\eta$ is of the form,
    \begin{align}
        \mathcal{C}_a (Q)=\{X\in \Omega_\eta: |X-Q|\leq a \ \text{dist}(X,\partial\Omega_\eta)\}.
    \end{align}

\end{definition}
\begin{definition}
    The non-tangential maximal function at $Q$ relative to $\mathcal{C}_a$ is
    \begin{align}
        Nu(Q)=\text{sup}\{|u(X)|:X\in \mathcal{C}_a(Q)\}.
    \end{align}
\end{definition}
We will omit the dependence of the cone on $a$ below.
\begin{definition}\label{Lp solvability} ($L^p$ solvability):
    The Dirichlet problem with data in $L^p(\partial\Omega_\eta, d\sigma)$ is solvable for $L$ if the solution for $Lu=0$ in $\Omega_\eta$, $u=f$ on $\partial\Omega$, for continuous boundary data $f$, satisfies the estimate
    \begin{align}
        ||N(u)||_{L^{p}(\partial\Omega_\eta,d\sigma)}\leq C ||f||_{L^{p}(\partial\Omega_\eta,d\sigma)}.
    \end{align}
\end{definition}

\begin{definition}\label{BMO}(BMO solvability):
    The Dirichlet problem is BMO-solvable for $L$ in the domain $\Omega_\eta$ if for all continuous $f$ with compact support on $\partial\Omega_\eta$, the solution $u$ of the classical Dirichlet problem with data $f$ satisfies the Carleson measure estimate
    \begin{align}
        \sup\limits_{x\in\partial\Omega_\eta,0<r<\infty} \frac{1}{\sigma(\Delta(x,r))}\int_{\Omega_\eta\cap B(x,r)} |\nabla u(y)|^2 \delta(y) dy\leq C ||f||^{2}_{BMO(\partial\Omega_\eta)}.
    \end{align}
\end{definition}
 For earlier work in this direction, we refer the reader to \cite{DKP11,Z18,HL18}.

 \begin{definition}
 ($A_{\infty}$, weak$A_\infty$): The elliptic measure associated with $L$ in $\Omega$ is of class $A_\infty$ with respect to the surface measure $\sigma=\mathcal{H}^{n}|_{\partial\Omega}$, which we denote by $\omega_L\in A_{\infty}(\sigma)$, if there exist $C_0 >1$ and $0<\theta<\infty$ such that for any surface ball $\Delta(q,r)=B(q,r)\cap \partial\Omega$, with $q\in \partial\Omega$ and $0<r<\text{diam}(\Omega)$, any surface ball $\Delta'=B'\cap \partial\Omega$ centered at $\partial\Omega$ with $B'\subset B(q,r)$, and any Borel set $F\subset\Delta'$, the elliptic measure with pole at the corkscrew point $A(q,r)$ relative to the surface ball $\Delta(q,r)$ satisfies
    \begin{align}
        \frac{\omega^{A(q,r)}_{L}(F)}{\omega^{A(q,r)}_{L}(\Delta')}\leq C_0\Bigg( \frac{\sigma(F)}{\sigma(\Delta')} \Bigg)^{\theta}.
    \end{align}

    On the other hand, if we have the following condition restricted to surface balls $2B'\subset B(q,r)$, and any Borel set $F\subset \Delta$
    \begin{align}
        \frac{\omega^{A(q,r)}_{L}(F)}{\omega^{A(q,r)}_{L}(2\Delta')}\leq C_0\Bigg( \frac{\sigma(F)}{\sigma(\Delta')} \Bigg)^{\theta},
    \end{align}
    then we say that the elliptic measure associated with $L$ in $\Omega$ , is of class weak $A_\infty$ with respect to the surfacem measure $\sigma$, and we denote this by $\omega_{L}\in \text{weak}A_\infty$.

    We further note that the weak $A_\infty$ condition implies the absolute continuity of $\omega^{A(q,r)}_{L}$ with respect to the surface measure $\sigma$ restricted to the surface ball $\Delta(q,r)$, and the following reverse H\"older estimate on the Radon-Nikodym derivative $k:=d \omega^{A(q,r)}_{L}/d\sigma$ for some $m>1$.
\quad  \begin{align} \left( \frac{1}{\sigma(\Delta)}\int_{\Delta} k^q \, d\sigma \right)^{1/m} \lesssim \frac{1}{\sigma(\Delta)}\int_{2\Delta} k \, d\sigma \approx \frac{\mu(2\Delta)}{\sigma(\Delta)}, \quad \forall \Delta' = B' \cap E, \text{ with } 2B \subseteq B(q,r).\end{align}
\end{definition}
\bigskip 
Given any arbitrary small $\eta$, we define the `ample' subdomains $B(0,1)\subset \R^{n+1}$ below. Here $2B$ is the euclidean ball with the same center as that of $B$ but with twice the radius.

\begin{definition}[Ample sawtooth domains]
For any $\eta$ small enough, a domain $\Omega_\eta\subset B(0,1)$ is an `ample' sawtooth domain with parameter $\eta$ if for each cube $Q_i\in \mathbb{D}_1$,   there exists a collection of disjoint sub-cubes $\{Q_{i,j} \subset Q_i \}$ with  $\sum_j \sigma(Q_{i,j})\leq \eta\sigma(Q_i)$ and we have, $\Omega_\eta= B(0,1)\setminus \sqcup_{i,j} T_{Q_{i,j}}$, where $T_{Q_{i,j}}$ is the Carleson box relative to the cube $Q_{i}$.   

\end{definition}

We use the terminology `ample' because these sawtooth domains, touch $\partial\Omega$ everywhere except for an $\eta$ fraction of the boundary.

We also give the definition of corkscrew domains.
\begin{definition}
    (Corkscrew domain):  An open set $\Omega\subset \R^{n+1}$ satisfies the interior corkscrew condition if for some uniform constant $c$ with $0<c<1$, and for every surface ball $\Delta:=\Delta(x,r)$ with $x\in \partial \Omega$ and $0<r<\text{diam}(\partial\Omega)$, there is a ball $B(A_{\Omega}(x,r),cr)\subset \Omega\cap B(x,r)$. The point $A_{\Omega}(x,r)\in \Omega$
 is called an interior corkscrew point relative to $\Delta$.
 \end{definition}

\begin{remark}
    We note that due to the background assumption \cref{large} on the drift, we have been able to use the Harnack inequality in the usual manner, throughout the domain. For reference, the form of the constant in the Harnack inequality in Theorem 8.20 in \cite{gt77}. Thus also, Harnack inequality will also be used throughout in subdomains as well.
\end{remark}

Now we state the main theorem in this paper.

\begin{theorem}\label{mainthm}
Let $L$ be the operator defined in \cref{operator} in $B(0,1)$. For any $\eta>0 $ arbitrarily small, the following holds: there exists an `ample' sawtooth domain $\Omega_{\eta}\subset B(0,1)$ with parameter $\eta$, so that assuming the Dirichlet problem for $L$ restricted to $\Omega_\eta$ is BMO solvable in $\Omega_\eta$, the elliptic measure $\omega_L$ belongs to weak-$A_\infty$ in the domain $\Omega_\eta$. 
    
\end{theorem}

\begin{remark}
    
The weak $A_\infty$ condition implies the solvability of the $L^p$ Dirichlet problem for some $p> 1$, for our operator $L$\footnote{This is in fact true for any open set $\Omega$ with Ahlfors-David regular boundary as mentioned in \cite{HL18}.}.See for example, Remark 2.22 in \cite{Dirichlet} which also considers this singular drift. The proof of this fact in our setting is essentially identical to the proof of Proposition 4.5 in \cite{HL18}, using the Boundary H\"older regularity obtained in our setting in the course of the proof of Claim 2 here, as a consequence of the Bourgain type estimate obtained of Claim 1. Thus essentially there is no distinction between the proof for the operator with the singular drift in our setting, versus the setting of an elliptic operator with no drift.
\end{remark}

\begin{remark}Note that the Boundary H\"older regularity inequality stated in Definition 9 contains a supremum over the ball $B(q,2r)$, $\sup\limits_{B(q,2r)\cap \Omega} u$, whereas in Corollary 2.4 of \cite{HL18} we have an average over the ball $B(q,2r)$. In our setting of the unit ball, and the fact that we can use Harnack's inequality because of the background assumption of \cref{large}, we also gain the Boundary Harnack principle (see for example, 4.4 of \cite{JK82}), thus finding an upper bound on $\sup\limits_{B(q,2r)\cap \Omega} u$ in terms of the value of $u$ at the corresponding corkscrew point, which is enough to conclude the proof of Proposition 4.5 of \cite{HL18}, in our unit ball setting.
\end{remark}

\bigskip

Here $\sigma$ is the surface measure on $\partial\Omega$, while $B(x,r)$ and $\Delta(x,r):=B(x,r)\cap\partial\Omega$ denote respectively the Euclidean ball in $\R^{n+1}$ and the surface ball on $\partial\Omega$ with center $x$ and radius $r$.

We recall the following result, as stated in \cite{HL18} Lemma 3.2, and proved in \cite{BL04}, which gives a sufficient condition for the elliptic measure to satisfy the weak-$A_\infty$ condition with respect to the surface measure. We use the same notation as in \cite{HL18}.

Given the point $x\in \Omega$, let $\hat{x}$ be a ``touching point" for the ball $B(x,\delta(x))$, i.e., $|x-\hat{x}|=\delta(x)$. We set, 
\begin{align}\label{delta}
    \Delta_{x}:=\Delta(\hat{x},10\delta(x)).
\end{align}

\begin{lemma}\label{lemma1}
    For any domain $\Omega$ with Ahlfors-David regular boundary, suppose that there are constants $c_0,\theta\in (0,1)$ such that for each $x\in\Omega$, with $\delta(x)<\text{diam}(\partial\Omega)$, and for every Borel set $F\subset \Delta_x$, we have,
    \begin{align}\label{eq14}
        \sigma(F)\geq (1-\theta)\sigma(\Delta_x) \implies \omega^{x}_L(F)\geq c_0.
    \end{align}
    Then we have $\omega^{y}\in \text{weak-}A_{\infty}(\Delta)$, where $\Delta=B\cap \partial\Omega$, for any ball $B=B(x,r)$, with $x\in\partial\Omega$ and $0< r<\text{diam}(\partial\Omega)$, and for all $y\in \Omega\setminus 4B$. The parameters in the weak-$A_\infty$ condition depend only on $n$, the ADR constants, $\eta,c_0$ and the ellipticity parameter $\lambda$ of the divergence form operator $L$.
\end{lemma}

Given the hypothesis of BMO solvability, and the Carleson measure estimate on the drift \cref{DKP2}, it will suffice to verify the hypothesis of Lemma 3, for the elliptic measure $\omega_L$. 

Given $\eta$, it will be enough to prove the two claims of Section 3 of \cite{HL18} in this setting, within the domain $\Omega_{\eta}$. We state them here.

Let $x\in \Omega_{\eta}$, and $\delta_{\Omega_{\eta}}(x)<\text{diam}(\partial\Omega_{\eta})$, and set $r:=\delta_{\Omega_{\eta}}(x)$. As before, choose $\hat{x}\in \partial\Omega_{\eta}$ so that, $|x-\hat{x}|=r$, and let $a$ be a sufficiently small number, depending only on $n$ and the ADR constants. Then we set,\footnote{Note that while we have defined the surface balls on the boundary $\partial\Omega_{\eta}$, we will continue to work with dyadic cubes in $S(0,1)$. }
\begin{align}\label{B_x}
    B_x := B(\hat{x},10r),\Delta_x =B_x \cap \partial\Omega_{\eta}. 
\end{align}
Then we have the following two claims, which are modifications of the ones in \cite{HL18} .
\begin{itemize}
\item[\textbf{Claim 1:}] Given $x\in \Omega_\eta$, for $a\leq 1/100$ a small enough uniform constant, there exists some $r_x$ with the property 
\begin{align}\label{impeq}
        \frac{\overline{c}\eps}{M}r\leq r_x\leq r,
 \end{align} 
so that there is a constant $\beta>0$ depending only on $n,\eps,c,M$, the ADR constants, and $\lambda$, and two balls $B_1 :=B(x_1,ar_x)\subset B_x, B_2:= B(x_2, ar_x)\subset B_{x}$, with $x_1, x_2\in \Delta_x$ such that $\text{dist}(x_1,x_2)> 5ar_x$, and 
\begin{align}\label{eq20}
    \omega_{L}^{x}(\Delta_1)\geq \beta , 
\end{align}
where $\Delta_1 := B_1 \cap \partial\Omega_{\eta}$, $\Delta_2 := B_2 \cap \partial\Omega_{\eta}$. 
We note that \cref{eq20} is effectively a Bourgain type estimate, given in \cref{Bourgain}.
\item[\textbf{Claim 2:}] Consider the function $u(x)$ that is a non-negative solution of $Lu=0$ in $\Omega_\eta$, vanishing continuously on $2\Delta_2$, with $||u||_{L^\infty}\leq 1$. Then for every $\gamma>0$ small enough, we have,
\begin{align}
    u(x)\leq C(\gamma)\Big(\frac{1}{\sigma(\Delta_{x})}\int_{B_x \cap\Omega} |\nabla u(y)|^2 \delta(y) dy\Big)^{\frac{1}{2}} +C \cdot \gamma^{\alpha}
\end{align}
\end{itemize}
for some exponent $\alpha$ that will be determined in the proof of this Claim in Section 2.3. Note that in effect the point $x_2$ plays the role that $\hat{x}$ plays the proof of \cite{HL18}. This is a technicality, since whenever the Whitney box $U_{Q_x}$ containing $x$, shares a boundary with $\partial\Omega_\eta$, we consider the points $x_1,x_2$, to lie on a lateral surface of $U_{Q_x}$ sharing the boundary with $\partial\Omega_\eta$, even though $\hat{x}$ may lie on a completely different boundary cube `below' $U_{Q_x}$, in $\partial\Omega_\eta$.

One is referred to earlier works that use versions of extrapolation of Carleson measures, which was originally constructed in \cite{LM95}, and versions of which were used in \cite{AHLT01,AHMTT02,HM12} , and which originally goes back to the Corona construction of \cite{Ca62, CG75}.

Our method of proof for both of these claims extends the method employed in \cite{HL18}. We use a stopping time argument, and work with a finite family of Whitney cubes, where each Whitney cube of each family fail the \cref{avg} condition. We have countably many such `bad cubes' to deal with, and we will need finitely many stopping time regimes. Furthermore, keeping in mind the background hypothesis of \cref{large}, we will be able to use the Harnack inequality and also finally use the Bourgain estimate within each of these bad Whitney cubes. These are explained in Section 2.

The main novelty in this paper is in dealing with the lack of Bourgain estimate of the elliptic measure in general, and thus we have to carefully extricate the `bad cubes' that fail \cref{avg} so that in the complement we have a collection of `good cubes' for which the Bourgain estimate does hold. We consider a nested sequence of domains on which we repeatedly use the Markov chain argument for the harmonic measure, and the main new ingredient is that we have to work in truncated balls of sufficiently small radius centered on the bottom `face' of the bad cubes, so that when restricted to such a truncated ball, we have the small constant regime for the drift and thus the Bourgain estimate for the harmonic measure for this truncated ball. After this we use a standard argument to relate the harmonic measure for this truncated ball, to the harmonic measure of the larger subdomain

As far as we are aware, this is the first time such a probabilistic extrapolation argument has been used where we have the dichotomy of `good' Whitney cubes with the pointwise smallness assumption on the drift and Bourgain estimate, and a set of bad Whitney cubes that fail the estimate of \cref{avg} , whose total volume is crudely bounded using the Carleson measure estimate. Subsequently, we need a modification to the Markov chain argument within each of the bad cubes.

This issue does not arise in the problem considered in \cite{FKP91, HM12}, where one assumes $L^p$ solvability of the Dirichlet problem for a divergence form elliptic operator, and with a Carleson measure assumption on the difference of the coefficient matrix, deduces the $L^p$ solvability of the Dirichlet problem for a perturbed divergence form elliptic operator. In such a problem, we uniformly have the Bourgain estimate for the perturbed operator, whereas for the drift term under consideration here, we only have the Bourgain estimate with the smallness assumption of \cref{smalll} below.\footnote{Further, we also lack uniform pointwise estimates on the Green's function for \cref{operator} in general, except again in the case of the pointwise small constant assumption on the drift term as in \cref{smalll}, whereas for the divergence form operator in the absence of the drift term, we have uniform pointwise Green's function estimates. }

Our argument after the conclusion of the proof of Claim 1 and Claim 2, is identical to that in \cite{HL18}, which we reproduce for the sake of completeness.


We state the following in a domain where the smallness assumption holds for each Whitney cube of the domain. For this result, see Chapter I on \cite{HL01} .
\begin{theorem}[Bourgain estimate (\cite{Bo87}):]\label{bourgain}
     Given any domain $\Omega\subset \R^{n+1}$, suppose we have the following pointwise smallness assumption on the drift for any $x\in \Omega$,
    \begin{align}\label{smalll}
        |\B(x)|\leq \frac{\eps}{\delta_{\Omega}(x)},
    \end{align}
   for some sufficiently small $\eps$. Then there exists a uniform constant $c>0$ so that for any $x\in \Omega$ where $|\hat{x}-x|=\delta_{\Omega}(x)=r$, $\Delta(x):= \Delta(\hat{x},10r)$ , we have 
   \begin{align}\label{Bourgain}
   \omega^{x}_{\Omega}(\Delta(x))\geq c
   \end{align}
   for the elliptic measure $\omega^{x}_{\Omega}$ corresponding to $L$.
\end{theorem}

We will use this Bourgain estimate in various subdomains constructed throughout the Markov chain argument; subdomains in which the drift satisfies the pointwise smallness assumption when the distance of the point is taken with respect to the boundary of the subdomain in consideration.

We note that the Bourgain estimate is equivalent to the Boundary H\"older continuity property, which we state below (see for example, Theorem 3.2 in Chapter II of \cite{HL01} for the implication that the Bourgain estimate implies Boundary H\"older continuity).

\begin{definition}\label{BHR}
   (Boundary H\"older regularity): There exist constants $C,\iota>0$ such that for $q\in \partial \Omega$ and $0<r<\text{diam}(\partial \Omega)$, and $u\geq 0$ with $Lu=0$ in $B(q,2r)\cap \Omega$, if $u$ vanishes continuously on $\Delta(q,2r):=B(q,2r)\cap \partial\Omega$, then,
    \begin{align}
        u(x)\leq C\Big( \frac{|x-q|}{r} \Big)^{\iota} \sup\limits_{B(q,2r)\cap \Omega} u \ \ \text{for any} \ x\in \Omega \cap B(q,r) . 
    \end{align}
\end{definition}

Finally we state that using the unique solvability of the continuous Dirichlet problem in subdomains where the drift term is bounded, it is well known that for the elliptic measure of a common segment of the boundary between two nested domains, with the same pole, gives a bigger value for the elliptic measure for the larger domain. This principle is used throughout the manuscript.





\bigskip



\section{Proofs of Claim 1 and Claim 2.}\label{Section2}
Before proving the two claims, we inductively outline the construction of the `ample' sawtooth domain with parameter $\eta$.

Let
\[
A = \{(i,k) : k \ge 1, i \in \mathfrak{J}_k, \text{ and } \int_P \sup_{y\in B(x, \delta(x)/2)} (|B|^2 \delta) dx > \epsilon \ell(P)^n \text{ for some } P \in \mathcal{W}' \text{ with } P \subset U_{Q_i^k}.\}
\]
That is, $A$ is the set of indices of Whitney cubes containing a ``bad'' cube. If $(i,k) \in A$, then for each $\hat{k} \le k$ there is a unique $\hat{i} \in \mathfrak{J}_{\hat{k}}$ with $Q_i^k \subseteq Q_{\hat{i}}^{\hat{k}}$; we then have that $(\hat{i}, \hat{k})$ either is or is not in $A$. We then let
\[
\mathfrak{K}_p = \{(i,k) \in A : p = \# \{ (\hat{i}, \hat{k}) \in A : Q_i^k \subseteq Q_{\hat{i}}^{\hat{k}} \}\}.
\]
We can then take the union over $\mathfrak{K}_p$ directly to define the $\F'_{p}, \F_{p}, \mathbb{G}_p$ as above.

We thus define,
\begin{align}
\F_{2} :=\sqcup_{(i,k)\in \mathfrak{K}_{2}}  T_{Q_{i}^{k}},
\end{align}

Inductively, we further define the families $\F_{p}$. By considering the union of the corresponding dyadic cubes in $S(0,1)$, we get the family,
\begin{align}
    \mathbb{G}_{p}:=\{ \cup Q_{m}^{k} :(m,k)\in \mathfrak{K}_{p}\}.
\end{align}

We also define analogously, using the index sets described above, for each $p\geq 1$, 
\begin{align}\label{eq25}
    \F'_{p} :=\sqcup_{(m,k)\in \mathfrak{K}_{p}}  S_{Q_{m}^{k}},  .
\end{align}

Now, because of the Carleson measure condition of \cref{DKP}, as well as \cref{large}, we will have at most 
\begin{align}\label{N_0}
N_0 :=\lceil (\lfloor2\sqrt{n+1}\rfloor)^{n} \mathscr{M}/\eps \eta \rceil
\end{align}

many such families of `bad ' Whitney cubes constructed above, till we reach $\F_{N_0}$ so that it has the property that 
\begin{equation}\label{impp}
\sum\limits_{Q_{i}\in \Gee_{N_0}} \sigma(Q_{i})\leq \eta \sigma(S(0,1)).\end{equation}

Note that we are using here the constant $\mathscr{M}$ from \cref{DKP}.

The value in \cref{N_0} comes about by noting that the families constructed above, $\mathbb{G}_k$, $k\geq 1$, are nested and the total surface measure of the union of the cubes in each of these $\mathbb{G}_k$ for fixed $k$, is bounded from below by $\eta \sigma(S(0,1))$. The corresponding cubes of $\mathcal{W}'$ , for the  fixed family with fixed $k$, thus contributes an amount of Carleson measure at least $\eps\eta/(\lfloor2\sqrt{n+1}\rfloor)^{n} \sigma(S(0,1))$. Thus the total number of families of  bad cubes one can construct in order to exhaust the total Carleson measure of $M\sigma(S(0,1))$, is at most the value of $N_0$ constructed above. 

Note of course, that we will in general have many Whitney cubes which satisfy \cref{avg} with positive Carleson measure, and so the union of the bad cubes (each of which by definition fail the \cref{avg} condition) will thus contribute an amount smaller than the amount given by $N_0$.

Next, we define the sawtooth domains, 
\begin{align}
    \Omega_{p}:= B(0,1)\setminus \F_{p}, \, \Lambda_p:= B(0,1)\setminus \F'_{p},\ \text{for} \ 1\leq p\leq N_{0}.
\end{align}

Note that, we have for each $1\leq p\leq N_0 -1$,
\begin{align}
     \Omega_p \subset \Lambda_p \subset \Omega_{p+1}.
\end{align}



\begin{figure}[h!]
\centering
\includegraphics[width=0.7\textwidth]{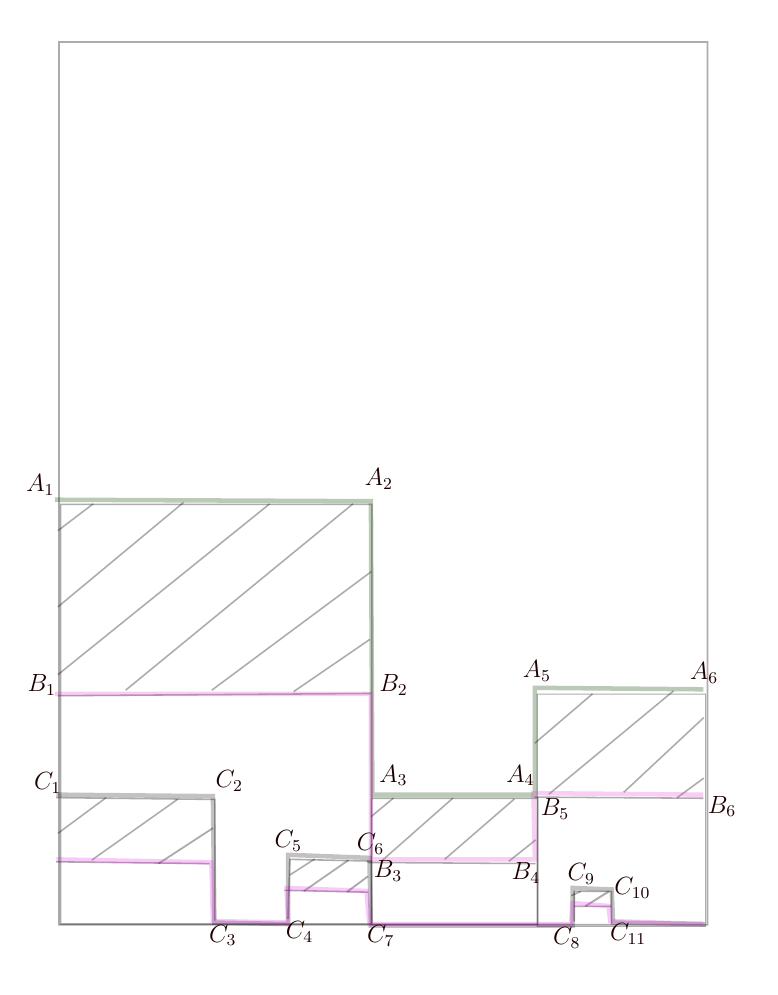}
\caption{For some $1\leq p\leq N_0 -1$, we show the intersection of the domains $\Omega_p, \Lambda_p$ to a given Carleson box $T_{Q^{1}_j}$ for a fixed $j\in \mathfrak{I}_1$. The domain $\Omega_p$ is delineated by the segments $A_1 A_2 A_3 A_4 A_5 A_6$ and lies above the union of these segments, while the domain $\Lambda_{p}$ is delineated by $B_1 B_2 B_3 B_4 B_5 B_6$ and lies above the union of these segments. The domain $\Omega_{p+1}$ is delineated by $C_1 C_2 C_3 C_4 C_5 C_6 C_7 C_8 C_9 C_{10} C_{11}$. The shaded regions denote the bad Whitney cubes which are the maximal cubes described in the construction in \cref{Section2}. Note in this figure that the points $C_6,B_3$ coincide, and the points $A_4,B_5$ coincide in the figure.}
\label{far field}
\end{figure}

\begin{proof}[Proof of Claim 1:]

We now consider the Markov argument for the elliptic measure corresponding to any pole $x\in \Omega_{N_0}$. First, consider the unique Whitney cube $U_{Q_x}\ni x$, with $U_{Q_x} \in \mathcal{W}(B(0,1))$.

It is enough to consider that $U_{Q_x}$ is a `bad' cube that contains at least one cube that fails the ASA condition. Otherwise if $U_{Q_x}$ is a good cube, the subsequent argument remains the same with minor modifications.

For the nested sequence of domains constructed above, consider the minimum value of $1\leq p_{x}\leq N_0$ so that $U_{x}\subset \Lambda_{p_{x}}$. Thus we have that $U_{Q_x}\subset\Lambda_{p_x}\subset \Omega_{p_{x}+1}$ and $U_{Q_{x}}\not\subset \Omega_{p_x}$.
 
 Then in the proof of Claim 1 as well as Claim 2, we will need to use the Markov process repeatedly, starting with the domain $\Omega_{p_{x}}$ , till the domain $\Omega_{N_0}$. 

 Let us note the Markov chain argument for the elliptic measure, for any fixed $1\leq p\leq N_0$,
\begin{equation}\label{eq27}
    \omega^{x}_{\Omega_{p}} (F_{p}) =\int_{y\in \partial\Lambda_{p-1}} \omega^{y}_{\Omega_{p}}(F_{p}) d\omega^{x}_{\Lambda_{p-1}}(y),  \ \omega^{x}_{\Lambda_{p}} (G_{p}) =\int_{y\in \partial\Omega_{p}} \omega^{y}_{\Lambda_{p}}(G_{p}) d\omega^{x}_{\Omega_{p}}(y)
\end{equation}
for any Borel subset $F_{p}\subset \partial\Omega_{p}$ and $G_{p}\subset \partial\Lambda_p$. 
 

\bigskip

Note that for the case of the unit ball $B(0,1)$ under consideration, it is enough to consider the Carleson boxes $\Tjj$ as these cover all of $B(0,1)$ except $B(0,1/2)$ . Whereas, if the pole of the harmonic measure lies inside $B(0,1/2)$ one can simply use Harnack inequality to compare the harmonic measure with poles anywhere in $\overline{B(0,1/2)}$.


Now we consider these two separate cases.
\begin{itemize}
    \item[Case 1.] First, consider the case where the Whitney box $U_x$ shares a boundary with $\partial \Omega_{N_0}$\footnote{Note that we can very well have a situation where we have a `bad' cube failing the ASA condition adjacent to $\Omega_\eta$.}. Let the length of this Whitney box be $l(U_{Q_x})=2^{-k_x}$.  In this case, consider a face of $\overline{U_{Q_x}}$ that intersects $\partial\Omega_{N_0}$, and call this $\delta_{Q_x,3}$\footnote{Recall that we denote by $\delta_{Q_x,1}, \delta_{Q_x,2}$ the two faces of $U_{Q_x}$ that are at constant radial distance from the origin. }.
    
    Without loss of generality, we may consider that $\hat{x}\in \delta_{Q_x,3}$, where $\hat{x}$ is the nearest point to $x$ on the boundary $\partial\Omega_{N_0}$. Otherwise, we will use Harnack's inequality within $U_{Q_x}$ to relate the harmonic measure with pole at $x$, to the harmonic measure with pole at a different point of $U_{Q_x}$ whose nearest point on $\partial\Omega_\eta$ lies on $\delta_{Q_x,3}$. Also, without loss of generality, it is enough to consider that $U_{Q_x}$ is a bad cube. Otherwise, if $U_{Q_x}$ is a `good' cube, then the argument simplifies as remarked at the end of the argument of Case 1. (See the remark following the statement of Claim 2).  

\bigskip

Consider,
\begin{align}\label{tau}
    \tau_{k_x}=\frac{\overline{c}\eps}{M}2^{-k_x},
\end{align}

and consider the ball $B(\hat{x},10a\tau_{k_x})$, for some uniform $a\leq 1/100$ in case $|x-\hat{x}|\geq a\tau_{k_x}$. Otherwise, if $|x-\hat{x}|< a\tau_{k_x}$, then we consider the ball $B(\hat{x},10r)$, where we recall that $r=|x-\hat{x}|$. For $U_{Q_x}\ni x$, and $\hat{x}\in \delta_{Q_x,3}$, we thus define,
  \begin{align}B(x)=\begin{cases}
   B(\hat{x}, 10a\tau_{k_x}), & \text{when}\ r\geq  a\tau_{k_x} \\
    B(\hat{x},10r), & \text{when}\ r< a\tau_{k_x}
\end{cases}\end{align}

Further consider the ball $B_{1}(x)$ with half the radius of the previous ball

    $B_1(x)=\begin{cases}
   B(\hat{x}, 5a\tau_{k_x}), & \text{when}\ r\geq  a\tau_{k_x} \\
    B(\hat{x},5r), & \text{when}\ r< a\tau_{k_x}
\end{cases}$

Note the difference in notation for the balls $B(x), B_1 (x)$, from that used in \cref{B_x}. Within the ball $B(x)$, we have the pointwise bound on the drift with the small constant, and thus from \cref{smalll,Bourgain} we have the Bourgain estimate. Note that with the choice of the constant $\tau_{k_x}$ in \cref{tau} is made precisely so that within the ball $B(x)$, the ambient bound on the drift with respect to the domain $B(0,1)$, yields the bound with the small constant within the ball $B(x)$, with respect to the domain $\Omega_{N_0}$, 
\begin{align}\label{eq40}
    |\B(y)|\lesssim \frac{M}{2^{-k_x}}\lesssim \frac{\eps}{\tau_{k_x}}\lesssim \frac{\eps}{\delta_{\Omega_{N_0}}(y)}, \ \ \forall y\in B(x)\cap \Omega_\eta.
\end{align}
    
The inequality above on the left comes from the ambient bound on the drift and the fact that $y\in U_{Q_x}$. The second inequality is a consequence of \cref{tau}, and the third follows immediately since $\delta_{\Omega_{N_0}(y)}\lesssim \tau_{k_x}$ for $y\in B(x)$.
\bigskip

Now define $t=\min(a\tau_{k_x},r)$.

We consider the two surface balls, $\Delta_1= B(x_1, at)\cap \delta_{Q_x,3} \subset \delta_{Q_x,3} \cap B_1 (x) , \Delta_2 = B(x_2, at )\cap \delta_{Q_x,3}\subset \delta_{Q_x,3} \cap B_1 (x) $ with,
 \begin{align}
   \text{dist}(\Delta_1,\Delta_2)\geq 2at.
\end{align}

 We have, the corkscrew point $A_{B(x)}(x_1,at)\in B(x_1,at)$ with respect to the domain $B(x)\cap \Omega_{\eta}$,  so that, the Bourgain estimate for the domain $B(x)\cap \Omega_\eta$, with pole at $A_{B(x)}(x_{1},at)$, gives us

\begin{align}\label{eq29}
\omega^{A_{B(x)}(x_{1},at)}_{B(x)\cap \Omega_\eta}(\Delta_1)\geq c.
\end{align}
Note that for the Bourgain estimate above, we needed to consider the surface balls $\Delta_1,\Delta_2$ to lie within $\delta_{Q_{x},3}\cap B_{1}(x)$ with the ball $B_1(x)$ of radius half of that of the ball $B(x)$ so that $\Delta_1,\Delta_2$  lie quantifiably away from the boundary of $B(x)$ and thus the Bourgain estimate with respect to $B(x)$, as in \cref{eq29}, holds  . Note that in the domain $B(x)$, if we denote the distance of a point $y\in B(x)$ to the boundary of the domain $\partial B(x)$ by $\delta_{B(x)}(y)$, then  for each $y\in B(x)$, similar to \cref{eq40}, we get,
\begin{align}\label{eq44}
    \frac{M}{2^{-k_x}}\leq \frac{\overline{c}\eps}{\tau_{k_x}}\leq \frac{\overline{c}\eps }{\delta_{B(x)}(y)}.
\end{align}
Thus we are in the small constant regime for the drift within the domain $B(x)$.

Then we also have, using the fact that the elliptic measure for the larger domain is bigger, for the common boundary segment,
\begin{align}\label{eq30}
\omega^{A_{B(x)}(x_{1},at)}_{\Omega_{\eta}}(\Delta_1)\geq \omega^{A_{B(x)}(x_{1},at)}_{B(x)\cap \Omega_\eta}(\Delta_1)\geq c. \end{align}

The right hand inequality above is the Bourgain estimate for the domain $B(x)$.

Further, we have, using Harnack's inequality, that there exists a uniform constant $c'=c'$ so that 
\begin{align}\label{eq31}
    \omega^{x}_{\Omega_{\eta}}(\Delta_1)\geq  c'\omega^{A_{B(x)}(x_{1},at)}_{\Omega_{\eta}}(\Delta_1).
\end{align}
Combining \cref{eq30,eq31}, we get 
\begin{align}\label{eq36}
    \omega^{x}_{\Omega_{\eta}}(\Delta_1)\geq c'c.
\end{align}
 Recall again that the stopping time argument was adopted and $N_0$ was chosen so that $\Omega_\eta =\Omega_{N_0}$. 

Note here that if the cube $U_{Q_x}$ was a `good' cube to begin with, then we would not need to consider the $\tau_{k_x}$ from \cref{tau}, since we will already have the bound
\
\begin{align}
    |\B(y)|\lesssim \frac{\eps}{2^{-k_x}}
\end{align}
for any $y\in U_{Q_x}$, and the radius of the surface balls $\Delta_1,\Delta_2$ can be simply taken to be comparable to $r=|x-\hat{x}|$, where as before, one can consider that $\hat{x}$ , lies on $\delta_{Q_x,3}$.
\bigskip 

\item[Case 2.] Next consider the case where the Whitney box $U_{Q_x}$ does not share a boundary with $\partial\Omega_{N_0}$. In this case, write in the coordinates of the unit ball, $x=(t_{x}, \theta_1,\dots,\theta_{n-1})$. Then, we consider the two points $x_1,x_2 \in S(0,t_x)\cap U_{Q_x}$ where $S(0,t_x)$ is the sphere of radius $t_x$ with center $0$, and with $\text{dist}(x_1,x_2)=5ar$. Recall again that $r=\partial_{\Omega_{\eta}}(x)=|x-\hat{x}|$, and so the length of the Whitney cube is comparable to $r$ and thus clearly a pair of such points $x_1,x_2$ exist. The points $x_1,x_2$ need not each be equidistant from $x\in U_{Q_x}$. In this case, using Harnack's inequality, we can relate the harmonic measures for subsets of $\partial\Omega_\eta$, with poles belonging anywhere in $U_{Q_x}$, by uniform constants.

Note that if $x'_1, x'_2$ are the radial projections of $x_1,x_2$ on the surface $\partial\Omega_{\eta}$, then we have by definition $x'_1,x'_2 \in \overline{T_{Q_x}}$ as well as $\text{dist}(x'_1,x'_2)\geq 5ar$. In this case, without loss of generality we  consider the surface ball $\Delta(x_1,ar)\subset S(0,t_x)\cap U_{Q_x} $, and consider the truncated sector $C_x:=(\{r:t_x \leq r\leq 1\}\times B') \cap \Omega_\eta$ where $B'$ is the radial projection of $\Delta(x_1,ar)$ on $S(0,1)$. Here, the surface ball $\Delta(x_1,ar)$ lies on the sphere $S(0,t_x)$, restricted to the Whitney cube $U_{Q_x}$. This truncated sector lies entirely within $\Omega_\eta$. We will perform the Markov chain argument within the truncated sector $C_x$. We write for each $p_x\leq p\leq N_0$, that, $C_x \cap \partial\Omega_p =F_p$, and $C_x\cap \partial\Lambda_p =G_p$, and finally that 
\begin{align} 
\Delta_1 =F_{N_0}.
\end{align} 

Note then that, with $x'_1$ as described above, we have $x'_1\in \Delta_1$.


Consider a point $v$ on the `bottom' face of the Whitney cube $U_{Q_x}$, which is the center of the surface ball $\delta_{Q_x,2}\cap C_x$, and then by using Harnack inequality, we will have, some uniform constant $C$ depending on the Harnack inequality constants, that, 
\begin{equation}
    \omega^{x_1}_{\Omega_{\eta}}(\Delta_{1})\geq C \omega^{v}_{\Omega_{\eta}}(\Delta_1).
\end{equation}

By another use of Harnack inequality, we also have that, 
\begin{align}
   C \omega^{x}_{\Omega_{\eta}}(\Delta_1)\geq \omega^{x_1}_{\Omega_{\eta}}(\Delta_1), \\
    C\omega^{x}_{\Omega_{\eta}}(\Delta_1)\geq \omega^{x_2}_{\Omega_{\eta}}(\Delta_1).
\end{align}

Henceforth, in this section, we write $\omega(\Delta)$in place of $\omega_{\Omega_{\eta}}(\Delta)$ when $\eta$ is understood.

Now consider the point $e_1\in \Omega_\eta$ on the radial line joining $0,v$, so that $\delta_{B(0,1)}(e_1)=ar$ , if such a point $e_1$ exists. If on the other hand, no such point $e_1$ exists, and instead there exists a point $e_2\in \partial\Omega_\eta$ on the radial line joining $0,v$ , so that $\delta_{B(0,1)}(e_2)> ar$, then $e_2$ is also the center of the surface ball $F_{N_0}$ defined above. Then consider the corkscrew point $A_{\Omega_\eta}(e_2,F_{N_0})$. Then in this second instance, we use Harnack's inequality to get, with a constant depending only on $a$, that, 
\begin{align}\label{eq500}
    \omega_{\Omega_\eta}^{v}(F_{N_0})  \geq C(a)\omega_{\Omega_\eta}^{A_{\Omega_\eta}(e_2,F_{N_0})}(F_{N_0}).
\end{align}

Further, using the Bourgain estimate, in either case where $A_{\Omega_\eta}(e_2,F_{N_0})$ lies in a `good' cube or it lies in a `bad' cube (in this latter case with a worse constant as in the analysis of Case 1 and in subsequent analysis), we will get a constant $c$ so that 
\begin{align}
    \omega_{\Omega_\eta}^{A_{\Omega_\eta}(e_2,F_{N_0})}(F_{N_0})\geq c.
\end{align}
This finishes the proof of Claim 1 in this case, and we have attained the bound of \cref{eq20}.

In the first instance where the point $e_1$ does exist, consider the minimum value of $1\leq u_x\leq N_0$ so that $e_1\in \Lambda_{u_x}$. Similar to \cref{eq500}, we now also have 
\begin{align}
    \omega_{\Omega_\eta}^{v}(F_{N_0})  \geq C(a)\omega_{\Omega_\eta}^{e_1}(F_{N_0}).
\end{align}

Henceforth we relabel $u_x$ as $p_x$ itself. We repeatedly use \cref{eq27}, starting with the domain $\Omega_{p_x +1}$, and initially the pole at $e_1$.

Recall that $S_{Q}$ was defined as the bottom half of the Carleson box corresponding to the surface cube $Q$. 


Using Harnack inequality, we consider the two following cases.
\begin{itemize} \item If $\Omega_{p_x +1} \cap C_x$ is nonempty, we do the following: consider the radial projection of $x_2$ onto the surface $\partial\Omega_{p_x +1}$, which is the center $x_{p_{x}+1}$ of the surface ball $F_{p_x +1 } \subset \partial\Omega_{p_x +1 }$ \footnote{Note that for each $p_x\leq p\leq N_0$, these surface balls can have subsets that coincide with subsets of  $S(0,1)$. }. For each $p_x\leq p\leq N_0$, we write $f_p= \text{diam}(F_{p}),\text{as well as} \  g_p =\text{diam}( G_p)$. Within the domain $S_{Q_x}$, we consider the surface ball $F'_{p_x+1 }$ with center $x_{p_x +1 }$  and diameter $\frac{1}{100} diam (F_{p_{x }+1})$. One considers the point $A_{p_x +1}=A_{S_{Q_x}}(x_{p_x+1}, \frac{1}{2}F'_{p_x +1})$ which is the interior corkscrew point corresponding to the surface ball $\frac{1}{2}F_{p_x +1}\cap \Omega_{p_x +1}$ within the corkscrew domain $ \Omega_{p_{x}+1}$, so that Bourgain's estimate gives us the second inequality on the right below,
\begin{align}\label{eq499}
 \omega_{\Omega_{p_x +1}}^{A_{p_x +1}}(\frac{1}{2}F_{p_x +1})   \geq  \omega_{B(x_{p_{x}+1}, \frac{1}{2}\text{diam}(F_{p_x  }))\cap \Omega_{p_x +1}}^{A_{p_x +1}}(F'_{p_x +1})\geq c.
\end{align}

The first inequality above was again obtained using the fact that the elliptic measure for the larger domain is bigger, for the common boundary segment.

Note again that we have chosen the surface ball $F'_{p_x +1}$ with diameter small enough so that the Bourgain estimate holds for the surface ball $F'_{p_x +1}$, restricted to the domain $B(x_{p_{x}+1}, \frac{1}{2}\text{diam}(F_{p_x }))\cap \Omega_{p_x +1}$. 

Lastly, we have, using Harnack's inequality, a constant $C'$ so that, 
\begin{align}\label{eq42}
   \omega_{\Omega_{p_x +1 }}^{e_1}(\frac{1}{2}F_{p_x })\geq C'\omega_{\Omega_{p_x +1 }}^{A_{p_x +1}}(\frac{1}{2}F_{p_x })\geq C' c=C'_1.
\end{align}

\item In case $ \Omega_{p_x +1} \cap C_x$ is empty, then we apply the argument for the domain $\Lambda_{p_x +1}$ as below.
\end{itemize}
Next, consider the Markov process for the domain $\Lambda_{p_x +1}$, and the Markov chain process,
\begin{align}\label{eq43}
   \omega^{e_1}_{\Lambda_{p_{x}+1}} (G_{p_{x}+1}) =\int_{y\in \partial\Omega_{p_{x}+1}} \omega^{y}_{\Lambda_{p_{x}+1}}(G_{p_{x}+1}) d\omega^{e_1}_{\Omega_{p_{x}+1}}(y) \geq \int_{y\in F_{p_{x}+1}} \omega^{y}_{\Lambda_{p_{x}+1}}(G_{p_x +1}) d\omega^{e_1}_{\Omega_{p_{x}+1}}(y).
\end{align}

Note that, $F_{p_x +1}$ can contain subsets that coincide with subsets of $S(0,1)$ or with `lateral' parts of Carleson boxes that constitute the union of $\mathbb{F}_{p_x +1}$, or on the `top' surfaces of Carleson boxes that constitute the union of $\mathbb{F}_{p_x +1}$.

\begin{itemize}
 \item  We estimate $\omega^{y}_{\Lambda_{p_x +1}}(G_{p_x +1})$ for any $y\in \delta_{Q,i}\subset F_{p_x +1} $ for $i\geq 3$ belonging to any `lateral' part of the top Whitney box $U_Q\subset T_Q$.  Let $l(Q)=2^{-k_Q}$. Further we recall the corresponding value from \cref{eq5}, of $\tau_k=\frac{\overline{c}\eps}{M}2^{-k_Q}$.

\begin{figure}
\centering
\includegraphics[width=0.28\textwidth]{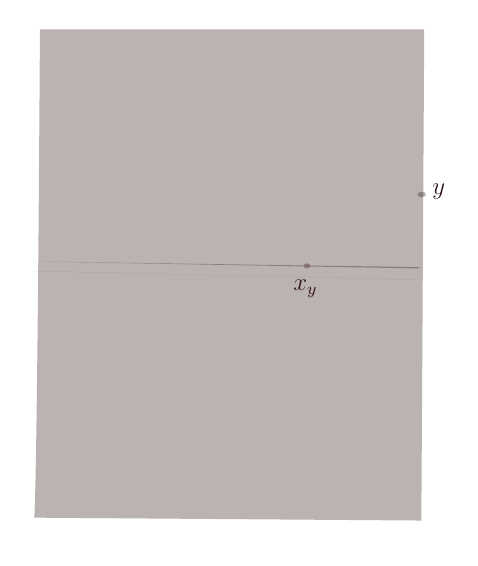}
\caption{The Carleson box $T_Q$ and the corresponding top half Whitney box $U_Q$, both corresponding to the surface cube $Q\in \mathbb{D}(S(0,1))$, are shown in the picture. The point $y\in \delta_{Q_i}$ for some $i\geq 3$ and the corresponding point $x_y \in \delta_{Q,2}$ on the bottom surface of $U_Q$ are shown in the figure. We have $\text{dist}(y,x_y)=2\delta(y)$, where $\delta(y)=\text{dist}(y,\delta_{Q,2})$.}
\label{fig 2}
\end{figure}

For each $y\in \delta_{Q,i}$ for $i\geq 3$, then we consider a ball $\Delta(x_y,a\delta(y)\frac{\overline{c}\eps}{M}):=\overline{B(x_y,a\delta(y)\frac{\overline{c}\eps}{M})}\cap \delta_{Q,2} \subset G_{p_x +1}$ where $\delta(y)=\text{dist}(y,\delta_{Q,2})$, is the distance of the point $y$ to the `bottom' surface $\delta_{Q,2}$ of the Whitney box $U_Q$. Further, we choose the point $x_y\in \delta_{Q,2} $, and require 
\begin{align}\label{eq53}
\text{dist}(y,x_y)=2\delta(y).
\end{align}

Consider the corkscrew point $A_{\Lambda_{p_x +1}}(x_y, a\delta(y)\frac{\overline{c}\eps}{M})$ corresponding to $B(x_y, a\delta(y)\frac{\overline{c}\eps}{M}) \cap \Lambda_{p_{x}+1}$ with respect to the domain $\Lambda_{p_{x}+1}$. In this case, using Harnack's inequality by connecting the points $y$ and $A_{\Lambda_{p_x +1}}(x_y, a\delta(y)\frac{\overline{c}\eps}{M})$ by a Harnack chain consisting of $\sim \log(\frac{M}{\overline{c}\eps})$  balls, we get, 
\begin{align}\label{eq444}
    \omega^{y}_{\Lambda_{p_x +1}}(G_{p_x +1})\geq \omega^{y}_{\Lambda_{p_x +1}}(\Delta(x_y, a\delta(y)\frac{\overline{c}\eps}{M})\geq c(M,\eps) \omega^{A_{\Lambda_{p_x +1}}(x_y, a\delta(y)\frac{\overline{c}\eps}{M})}_{\Lambda_{p_x +1}}(\Delta(x_y, a\delta(y)\frac{\overline{c}\eps}{M})
\end{align}
We note that by construction, the corkscrew point belongs to the ball $B(x_y, a\frac{\overline{c}\eps}{M}2^{-k_Q})$ where we can use the estimates for the small constant regime for the drift, with respect to the domain $\Lambda_{p_x +1}$, for the following reason. Since the bound on the drift in $B(0,1)$ is given by $ M/\delta_{B(0,1)}(z)$, this can be bounded from above as,
\begin{align}
    \frac{M}{\delta_{B(0,1)}(z)}\leq \frac{\overline{c}\eps}{\delta_{ B(x_y, a\frac{\overline{c}\eps}{M}2^{-k_Q})\cap \Lambda_{p_x +1}}(z)}, \text{for } z\in B(x_y, a\frac{\overline{c}\eps}{M}2^{-k_Q})\cap \Lambda_{p_x +1}.
\end{align}

Thus finally, using the Bourgain estimate for the domain $ B(x_y, a\frac{\overline{c}\eps}{M}2^{-k_Q})\cap \Lambda_{p_x +1}(z)$ and the maximum principle, which gives us that the harmonic measure for a segment of a subset of the common boundary of two nested domains, is greater for the larger of the two nested domains, we get, 
\begin{align}\label{eq455}
    \omega^{A_{\Lambda_{p_x +1}}(x_y, a\delta(y)\frac{\overline{c}\eps}{M})}_{\Lambda_{p_x +1}}(\Delta(x_y, a\delta(y)\frac{\overline{c}\eps}{M}))\geq \omega^{A_{\Lambda_{p_x +1}}(x_y, a\delta(y)\frac{\overline{c}\eps}{M})}_{B(x_y, a\frac{\overline{c}\eps}{M} 2^{-k_Q})\cap\Lambda_{p_x +1}}(\Delta(x_y, a\delta(y)\frac{\overline{c}\eps}{M})   \geq c.
\end{align}
Thus we have for each $y\in \delta_{Q,i}$ for $i\geq 3$, using the above equations,
\begin{align}\label{eq47}
    \omega^{y}_{\Lambda_{p_x +1}}(G_{p_x +1})\geq c(M,\eps)c.
\end{align}

 \item For any $y$ belonging to the `top' surface $\delta_{Q,1}$ of the Whitney box $U_Q$, we consider the point $x_y$ on the `bottom' surface $\delta_{Q,2}$ of $U_Q$, and on the same radial line connecting $0,y$. Thus we have using the Harnack inequality, 
\begin{align}\label{eq57}
    \omega^{y}_{\Lambda_{p_x +1}}(G_{p_x +1})\geq \omega^{y}_{\Lambda_{p_x +1}}(\Delta(x_y, a\delta(y)\frac{\overline{c}\eps}{M})\geq c(M,\eps) \omega^{A_{\Lambda_{p_x +1}}(x_y, a\delta(y))}_{\Lambda_{p_x +1}}(\Delta(x_y, a\delta(y)\frac{\overline{c}\eps}{M})
\end{align}
As in \cref{eq455}, we further get that, 
\begin{align}\label{eq49}
    \omega^{A_{\Lambda_{p_x +1}}(x_y, a\delta(y))}_{\Lambda_{p_x +1}}(\Delta(x_y, a\delta(y)\frac{\overline{c}\eps}{M})\geq \omega^{A_{\Lambda_{p_x +1}}(x_y, a\delta(y)\frac{\overline{c}\eps}{M})}_{B(x_y, a\delta(y)\frac{\overline{c}\eps}{M})\cap \Lambda_{p_x +1}}(\Delta(x_y, a\delta(y)\frac{\overline{c}\eps}{M})   \geq c.
\end{align}
Thus combining the above two estimates, we have, 
\begin{align}\label{eq50}
     \omega^{y}_{\Lambda_{p_x +1}}(G_{p_x +1})\geq c(M,\eps)c
\end{align}

\item For the parts of $F_{p_x +1}$ that coincide with $S(0,1)$, or parts of $F_{p_x +1}$ that coincide with the `lateral' surfaces of $S_Q$(the bottom half of the Carleson box $T_Q$ corresponding to $Q$) we have a sharp factor of $1$ for the elliptic measure. In these cases, we move on to the next step in the Markov chain argument. 

\end{itemize}

Thus finally, using \cref{eq43}, and combining the estimates from \cref{eq47,eq49}, and the crude factor of $1$ from the subsets of $F_{p_x +1}$ that belong to $S(0,1)$ or some lateral surface of $S_Q$ for some $Q$, we have, 
\begin{align}\label{eq51}
    \omega_{\Lambda_{p_x +1}}^{e_1}(G_{p_x +1})=\int_{y\in \partial\Omega_{p_x +1}} \omega^{y}_{\Omega_{p_x +2}}(G_{p_{x}+1} )d\omega^{e_1}_{\Lambda_{p_x +1}}(y)\\ \geq \int_{y\in F_{p_x +1}} \omega^{y}_{\Omega_{p_x +2}}(G_{p_{x}+1} )d\omega^{e_1}_{\Lambda_{p_x +1}}(y)\geq c_{1}(M,\eps)C'_1,
\end{align}

where in the last step, we used the estimate from \cref{eq42}.

From here, using the version of the left hand side of \cref{eq27} using the pole at $e_1$ in place of $x$, now to estimate $\omega^{e_1}_{\Omega_{p_x +2}}(F_{p_x +2})$, we consider the following separate cases, 
\begin{itemize}
    \item For $y\in \partial \Lambda_{p_x +1}$, suppose we have for some $Q\in \mathbb{D}$, that $y\in \partial S_Q $, with $y\in \delta'_{Q,i}$ with some $i\geq 3$, i.e. some `lateral' surface of $S_{Q}$ . We consider the two separate cases,
    \begin{itemize}
    \item If $y\in \partial T_{Q_1}$ for some $T_{Q_1}\in \mathbb{F}_{p_{x+2}}$ (recall the construction of $\mathbb{F}_{p_{x+2}}$ in the beginning of Section 2 ), then we pick up a crude factor of $1$ in this step. Then we will move further to the subsequent step in the Markov chain process.
    
    \item Otherwise, $y\notin \partial T_{Q_1}$ for any $T_{Q_1}\in \mathbb{F}_{p_{x+2}}$, and consider the point $x'_y\in \text{int}(S_Q)\cap \partial \Omega_{p_x +2} $ in the interior of $S_Q$, defining $\text{dist}(y,\partial\Omega_{p_x +2}):=\delta_1 (y)$, with the property that, 
   \begin{align}\label{eq62}
    2 \delta_1 (y)\geq \text{dist}(x'_y,\delta'_{Q,i})\geq \delta_1 (y)\end{align}
    
    for some $i\geq 3$. This point $x'_y $ may belong to either a `top' or `lateral' surface of a Carleson box belonging to $S(0,1)\setminus \Omega_{p_x +2}$ , or on $S(0,1)$ itself. The situation is analogous to the one depicted in \cref{fig 2} for the previous case. We thus consider the corkscrew point $A_{\Omega_{p_x +2}}(x'_y, a\delta_1 (y))$ corresponding to the surface ball $\Delta(x'_y, a\delta_1 (y)):=B(x'_y,a\delta_1 (y))\cap \Omega_{p_x +2}$ .Recall that $a\leq 1/100$ is a uniformly small constant.

    Note that, by this construction, for any $ \delta'_{Q,i}\ni y$, with any $i\geq 3$, we get a uniform bound on the length of the Harnack chain joining $y$ with $A_{\Omega_{p_x +2}}(x'_y, a\delta_1 (y))$, and thus the bounds on the Harnack inequality are uniform in all such $y$.

    \end{itemize}
 
Thus we have, 
\begin{align}\label{eq44}
    \omega^{y}_{\Omega_{p_x +2}}(F_{p_x +2})\geq \omega^{y}_{\Omega_{p_x +2}}(\Delta(x'_y, a\delta_1 (y))\geq c'(M,\eps) \omega^{A_{\Omega_{p_x +2}}(x'_y, a\delta(y))}_{\Omega_{p_x +2}}\Delta(x'_y, a\delta_1 (y)),
\end{align}

for some constant $c'(M,\eps)$.

 We can now use the Bourgain estimate for the domain $ B(x_y, a2^{-k_Q})\cap \Omega_{p_x +2}$ and the maximum principle, to get that, 
\begin{align}\label{eq45}
    \omega^{A_{\Omega_{p_x +2}}(x'_y, a\delta_1 (y))}_{\Omega_{p_x +2}}(\Delta(x'_y, a\delta_1(y))\geq \omega^{A_{\Omega_{p_x +2}}(x'_y, a\delta_1 (y))}_{B(x'_y, a\delta_1(y))\cap \Omega_{p_x +2}}(\Delta(x'_y, a\delta_1 (y))   \geq c.
\end{align}

Combining the above two estimates, we get, 
\begin{align}\label{eq54}
    \omega^{y}_{\Omega_{p_x +2}}(F_{p_x +2})\geq  c'(M,\eps)\omega^{A_{\Omega_{p_x +2}}(x'_y, a\delta(y))}_{B(x'_y, a\delta_1(y))\cap \Omega_{p_x +2}}(\Delta(x'_y, a\delta_1 (y))   \geq c'(M,\eps)c.
\end{align}

\item Now when $y\in \delta'_{Q,1}$, the top surface of $S_Q$, then we again consider the two cases, 
\begin{itemize}
\item If $y\in \partial T_{Q_1}$ for some $T_{Q_1}\in \mathbb{F}_{p_x +2}$, then we pick up a crude factor of $1$ in this step. Then we will move further to the subsequent step in the Markov chain process.

\item Otherwise, consider a point $x'_y \in \text{int} (S_Q) \cap \partial \Omega_{p_x +2}$, on the same radial line through the origin $0$ and $y$, with the property that,
\begin{align}\label{eq665}
\text{dist}(x'_y, \delta'_{Q,1})\geq \frac{1}{2}2^{-k_Q -1}
\end{align}

 then, again by arguments similar to \cref{eq49,eq50}, 
\begin{align}
     \omega^{y}_{\Omega_{p_x +2}}(F_{p_x +2})\geq \omega^{y}_{\Omega_{p_x +2}}(\Delta(x'_y, a\delta_1 (y))\geq c'(M,\eps) \omega^{A_{\Omega_{p_x +2}}(x'_y, a\delta_1 (y))}_{\Omega_{p_x +2}}(\Delta(x'_y, a\delta_1 (y))).
\end{align}

Note that, in this process, the Harnack chain joining $y$ with $A_{\Omega_{p_x +2}}(x'_y, a\delta_1 (y))$ has length that is uniform for any choice of $y\in \partial\delta'_{Q,1}$, by the construction of $A_{\Omega_{p_x +2}}(x'_y, a\delta_1 (y))$.

Then, we also have by the Bourgain estimate for the domain $ B(x'_y, a\delta_1 (y))\cap \Omega{p_{x}+2}$, the right hand estimate below,
\begin{align}\label{eq56}
    \omega^{A_{\Omega_{p_x +2}}(x'_y, a\delta_1 (y))}_{\Omega_{p_x +2}}(\Delta(x'_y, a\delta_1(y))\geq \omega^{A_{\Omega_{p_x +2}}(x'_y, a\delta_1 (y))}_{B(x'_y, a\delta_1(y))\cap\Omega_{p_x +2}}(\Delta(x'_y, a\delta_1 (y))   \geq c.
\end{align}
Again we have used the fact that the elliptic measure is greater for the larger domain, for the common segment of the boundary of the two domains for the first inequality above.

And thus we have that, 
\begin{align}\label{eq57}
     \omega^{y}_{\Omega_{p_x +2}}(F_{p_x +2})\geq c'(M,\eps)c:=c_2 (M,\eps).
\end{align}

\end{itemize}

\end{itemize}

Thus collecting all these estimates, and using the right hand side of \cref{eq27} , we will get using \cref{eq54,eq57}, and the crude estimate of $1$ from the other two cases, that, 
\begin{multline}
    \omega^{e_1}_{\Omega_{p_x +2}}(F_{p_x +2})=\int_{y\in \partial\Lambda_{p_x +1}} \omega^{y}_{\Omega_{p_x +2}}(F_{p_{x}+2} )d\omega^{e_1}_{\Lambda_{p_x +1}(y)}\geq \int_{y\in G_{p_x +1}}\omega^{y}_{\Omega_{p_x +2}}(F_{p_{x}+2} )d\omega^{e_1}_{\Lambda_{p_x +1}}(y)\\ 
  \geq c_2(M,\eps) c_1(M,\eps) C'_1,
\end{multline}
where the last step follows from the estimate in \cref{eq51}.

Note that $\int_{y\in G_{p_x +1}}\omega^{y}_{\Omega_{p_x +2}}(F_{p_{x}+2} )d\omega^{e_1}_{\Lambda_{p_x +1}}(y)$ is an upper bound for the actual contribution we are considering for a set of $y$ values in a proper subset of $G_{p_x +1}$.





Then the argument continues analogously, by alternating between the two sets of domains $\Omega_p,\Lambda_p$ till we reach the domain $\Omega_{N_0}=\Omega_\eta$.

Along with the above Markov chain estimate, we separately note that starting with the choice of the point $e_1$ made so that $\delta_{B(0,1)}(e_1)=ar$, and further the choices of the points $x_y, x'_y$ made in \cref{eq53,eq62,eq57,eq665} at each of the $N_0 -p_x$ steps, the choice of $\text{diam}(F'_{p_{x}+1})=\frac{1}{100}\text{diam}(F_{p_x +1})$, we ensure the requisite lower bound on the harmonic measure for $\Delta_1$. In other words, the harmonic measure does not `spill out' beyond the cone $C_x\cap \partial\Omega_\eta=\Delta_1$. By repeatedly using \cref{eq53,eq62,eq57,eq665}, we get that this harmonic measure is restricted within $C_x \cap \partial\Omega_\eta$ and we get the estimate of \cref{eq60} below. In fact, the choice of $e_1$ prior to \cref{eq500} was made precisely for this purpose.

\bigskip

In the process, for any point $x\in \Omega_{\eta}$ where $U_{Q_x}$ does not share a boundary with $\partial\Omega_\eta$, we will pick up a lower bound on the harmonic measure, of the form,
\begin{align}\label{eq59}
    \omega^{x}_{\Omega_{\eta}}(\Delta_1)\geq C\omega^{e_1}_{\Omega_\eta}(\Delta_1)\geq CC_1' \big( c_2(M,\eps) c_1(M,\eps)\big)^{N_0 -p_x}\geq C(a) (c_2(M,\eps) c_1(M,\eps))^{N_0}.
\end{align}
 Combining the two cases, from \cref{eq36,eq59}, we get a uniform $C_3 (\eta)$ dependent on $\eta$ so that for any $x\in \Omega_\eta$, we have, 
\begin{align}\label{eq60}
    \omega^{x}_{\Omega_\eta}(\Delta_1)\geq C_3(\eta) .
\end{align}

\end{itemize}
This concludes the proof of Claim 1.
\end{proof}

\bigskip
\begin{proof}[Proof of Claim 2]
    This follows by using the Markov property, Claim 1, and finally the corresponding argument for Claim 2 from \cite{HL18}.  We separate the two cases, as before, depending on whether $U_{Q_{x}}$ shares a boundary with $\Omega_\eta$ or not. We note that this result is essentially a consequence of the boundary H\"older regularity, which for the first case already holds, and which for the second case needs to be proved using the result of Claim 1.
    \begin{itemize}
        \item[Case 1.] Suppose that $U_{Q_x}$ shares a boundary with $\Omega_\eta$. In that case, we recall the construction from Claim 1. In this case, we use boundary H\"older continuity in the domain $B(x)$ (see for example, the Remark after the proof of Lemma 3.9 in Chapter I of \cite{HL01}.). 
We follow the argument of the proof of Claim 2 in \cite{HL18}. Recall the notation introduced in this case, for the proof of Claim 1. Without loss of generality, we consider that $x_2=0$, and the open horizontal cone $\Gamma$ of spherical aperture $\pi/100$ which intersects $U_{Q_x}$ and then the spherical cap inside $S:S^n\cap \Gamma$. Consider the corkscrew point $X:=A_{B(x)}(x_2, at)$ and $r'_x=|X-x_2|$. Then we get,
\begin{align}
    u(X)\lesssim \int_{S} u(r'_x\xi)d\mu(\xi) =\int_{S}\Big( u(r'_\xi)-u(\gamma r'_\xi) d\mu (\xi) +O(\gamma^{\alpha})\Big),
\end{align}

where for the last term on the right, we used the fact that the boundary H\"older regularity holds in $B(x)$ along with the fact that $||u||_{L^{\infty}(\Omega)}\leq 1$.
Then, we get,
\begin{align}
    |I|=\Big| \int_{S}\int_{\gamma r'_{x}}^{r'_{x}} \frac{\partial}{\partial t}(u(t\xi))dt d\mu(\xi)\Big| \leq (\gamma r'_x)^{-n}\int\int_{\Gamma\cap (B(0,r'_x)\setminus B(0,\gamma r'_x))} |\nabla u(y)| dY.
\end{align}

where we used polar coordinates in $(n+1)$ dimensions. Now we have,
\begin{multline}
   u(X)\leq c \gamma ^{-n -1/2}(r'_{x})^{-n/2}\Big(  
 \int\int_{B(x_2, at)\cap \Omega_{\eta}}|\nabla u(y)|^{2} \delta(y)dy\Big)^{1/2} +O(\gamma^\alpha)\\ \leq C(a_o,\eps,M) \gamma^{-n-1/2} \delta(x)^{-n/2}\Big(  
 \int\int_{B(x_2, at)\cap\Omega_{\eta}}|\nabla u(y)|^{2} \delta(y)dy\Big)^{1/2} +O(\gamma^\alpha)\\ \leq C(a_o,\eps,M) \gamma^{-n-1/2} \delta(x)^{-n/2}\Big(  
 \int\int_{B(x,10\delta(x))\cap\Omega_\eta}|\nabla u(y)|^{2} \delta(y)dy\Big)^{1/2}+O(\gamma^{\alpha}),
\end{multline}
where in the last two steps, we have reverted back to the point $x\in U_{Q_x}$ , using a factor dependent on $a,\eps,M$ to change from $r'_x$ to $\delta(x)$. 
Next we can use the Harnack inequality with a Harnack chain connecting $x, X$ of length dependent on $a,\eps,M$, so that we have,
\begin{multline}
    u(x)\leq C_1(a,\eps,M)u(X)\\ \leq C(a_o,\eps,M) \gamma^{-n-1/2} \delta(x)^{-n/2}\Big(  
 \int\int_{B(x,10\delta(x))\cap\Omega_\eta}|\nabla u(y)|^{2} \delta(y)dy\Big)^{1/2}+O(\gamma^{\alpha}).
\end{multline}
\item[Case 2.] In case $U_{Q_x}$ does not share a boundary with $\Omega_\eta$, we follow a slightly different approach below.

Consider any $q\in \Delta_2$ and any $0<r< \text{diam}(\Omega)$, and the domains $D(q,k,r):= B(q,10^{-k}r)\cap  \Omega_\eta $ for all integers $k\geq 0$, and consider the boundaries of these domains that intersect the interior of $\Omega_\eta$, i.e. $\partial_1 D(q,k,r):=S(q,10^{-k}r) \cap \text{int}(\Omega_\eta)$, as well as the part of the boundary of $D(q,k)$ coinciding with $\partial\Omega_\eta$, which we write $\partial_2 D(q,k,r):= \partial D(q,k,r)\cap \partial \Omega_\eta$. 

Consider any $k\geq 1$. For any point $y_k \in \partial_1 D(q,k,r)$, we note that Claim 1 is adapted immediately for the domain $D(q,k-1,r)$ (which contains $D(q,k,r)$), and further using Harnack inequality, we have, 
 \begin{align}\label{eq65} \omega_{D(q,k-1,r)}^{y_k}(\partial_2 D(q,k-1,r))\geq \beta,\end{align} 
for some uniform $\beta$ independent of $k$.

Further, we also have, by the unique solvability of the Dirichlet problem in these subdomains, that for each $k\geq 1$,
\begin{align}\label{eq66}
    u(y_k)=\int_{\partial_1 D(q,k-1,r) } u(y_{k-1})d\omega_{D(q,k-1,r)}^{y_k}(y_{k-1}) +\int_{\partial_2 D(q,k-1,r) } u(y_{k-1})d\omega_{D(q,k-1,r)}^{y_k}(y_{k-1}).
\end{align}
Note that the second term on the right is identically $0$, since $u=0$ on the $\partial_2 D(q,k,r)$ by hypothesis.

We consider that $||u||_{L^{\infty}}\leq 1$, by the hypothesis of Claim 2.

By definition of the elliptic measure, we have, 
\begin{align}\label{eq67}
    \int_{\partial_1 D(q,k-1,r) } d\omega_{D(q,k-1,r)}^{y_k}(y_{k-1}) +\int_{\partial_2 D(q,k-1,r) } d\omega_{D(q,k-1,r)}^{y_k}(y_{k-1})=1
\end{align}
Combining \cref{eq65,eq66,eq67}, noting again that $||u||_{L^{\infty}}\leq 1$ in $\Omega_\eta$, we get,
\begin{multline}
     u(y_k)=\int_{\partial_1 D(q,k-1,r) } u(y_{k-1})d\omega_{D(q,k-1,r)}^{y_k}(y_{k-1})\\ \leq \text{sup}_{\partial_{1}D(q,k,r)} \ u(y_{k-1})\int_{\partial_1 D(q,k-1,r) } d\omega_{D(q,k-1,r)}^{y_k}(y_{k-1})\leq (1-\beta)\text{sup}_{y_{k-1 }\in \partial_{1}D(q,k-1,r)} u(y_{k-1})
\end{multline}

For any $r_1$, choose $k$ so that $r_1\approx r 10^{-k}$, and then we get from above for $|y_k -q|\approx r_1$
\begin{align}
    u(y_k)\leq (1-\beta)^{k}\text{sup}_{y_1 \in \partial_1 D(q,1,r)}u(y_1)\leq (1-\beta)^{k}.
\end{align}
From here, we immediately get the 
\begin{align}
    u(y_k)\leq C \big(\frac{r_1}{r}\big)^{\alpha},
\end{align}
for some uniform $\alpha>0$, thus establishing the Boundary H\"older continuity. After this, the inequality of Claim 2 follows immediately as in the previous case.
    \end{itemize}
\end{proof}

After this, the argument follows exactly as in \cite{HL18}, using Lemma 3.9 exactly as stated in that paper. That argument in turn follows the argument in \cite{DKP11}. Note that the argument here remains identical in our situation with the operator with the drift term. For the sake of completeness, we reproduce the argument of \cite{HL18} here, to establish the hypothesis of Theorem 3. Let $F\subset \Delta_x$ be a Borel set satisfying the left hand condition \cref{eq14}. If $\theta$ is small enough, depending on $\eps,n,M, c$, the ADR constants, so that, 
\begin{align}
    \sigma(F_1)\geq (1-\sqrt{\theta})\sigma(\Delta_1),
\end{align}
where $F_1:=F\cap \Delta_1$. Set $A_1 :=\Delta_1 \setminus F_1$, and define, 
\begin{align}
    f:=\text{max}(0,1+\gamma \log \mathcal{M} (1_{A_1})), 
\end{align}
where $\gamma$ is a small constant to be chosen, and $\mathcal{M} $ is the Hardy-Littlewood maximal operator on $\partial\Omega$. We note that,
\begin{align}
    0\leq f\leq 1,\ \ \ ||f||_{BMO(\partial\Omega)} \leq C\gamma, \ \ 1_{A_1}\leq f.
\end{align}

Also note that if $z\in \partial\Omega \setminus 2B_1$, then
\begin{align}
    \mathcal{M}(1_{A_1})(z)\lesssim \frac{\sigma(A_1)}{\sigma(\Delta_1)}\lesssim \sqrt{\theta}.
\end{align}
Then if $\theta$ is chosen small enough depending on $\gamma$, then $1+\gamma\log(\mathcal{M})(1_{A_1})$ will be negative, hence $f=0$ on $\partial\Omega\setminus 2B_1$.

In order to work with continuous data as in the definition of BMO solvability, we state the following lemma, 
\begin{lemma}\label{lemma5}
    There exists a collection of continuous functions $\{f_s \}_{0<s<1}$ defined on $\partial\Omega$, with the following properties,
    \begin{enumerate}
        \item $0\leq f_s\leq 1$ for each $s$.
        \item $\text{supp}(f_s)\subset 3B_1 \cap \partial\Omega$.
        \item $1_{A_1}(z)\leq \text{lim inf}_{s\to 0} f_s (z)$, for $\omega^{X}-a.e\  z\in \partial\Omega$.
        \item $\text{sup}_{s}||f_s||_{BMO(\partial\Omega)}\leq C||f||_{BMO(\partial\Omega)}\lesssim \gamma$, where $C=C(n,ADR)$.
    \end{enumerate}
\end{lemma}

We do not reproduce the proof of this lemma here, and one is referred to \cite{HL18}.

Let $u_s$ be the solution of the Dirichlet problem for the equation $L u_s =0$ in $\Omega$, with continuous data $f_s$. Then for a small $\eps>0$, to be chosen momentarily, by Fatou's lemma, and Claim 2, along with Lemma 5, we have, 
\begin{align}\label{3.10}
    \omega_{L}^{x}(A_1)\leq \int_{\partial\Omega}\text{lim inf}_{s\to 0} f_s d\omega^{X} \leq \text{lim inf}_{s\to 0} u_s (X)\leq C_\eps \gamma +\overline{c}\eps^{\alpha},
\end{align}

Next set $A:=\Delta_x \setminus F$, and by the definition of $A,A_1$, along with Claim 1, and \cref{3.10}, we get that,
\begin{align}
    \omega_{L}^{x}(A)\leq \omega_{L}^{x}(\Delta_x \setminus\Delta_1)+\omega_{L}^{x}(A_1)\leq (\omega_{L}^{x}(\Delta_x)-\beta +C_\theta \gamma +\overline{c}\theta^{\alpha}).
\end{align}

We now choose first $\theta>0$ and then $\gamma>0$, so that, $C_\theta\gamma +\overline{c}\theta^{\alpha}<\beta/2$, to obtain, 
\begin{align}
    \omega_{L}^{x}(\Delta_x)-\omega_{L}^{x}(A)=\omega_{L}^{x}(F)\geq \frac{\beta}{2}.
\end{align}

Recall that $\beta$ depends on the ambient parameters of the problem as defined in the statement of Claim 1. We have thus proved the hypothesis of  Lemma 3 for this problem, and thus have established the weak-$A_\infty$ property.

\begin{remark}
    Note that we have presented the argument for the unit ball, for simplicity of notation. However, the methods can be extended in the expected manner to more general bounded Lipschitz domains with small constants, keeping track of the `top', `bottom' and `lateral' faces of the Whitney cubes throughout the Markov chain argument. This will be investigated in future.
\end{remark}


\begin{remark}
    Note that we have strictly only used the Carleson measure condition of \cref{DKP} for the cubes in $\mathbb{D}_1$, in constructing the sequence of subdomains. Our result is not a scale invariant estimate, and depends on the finite volume of the ball, and the result is the specified 'ample' subdomains.
\end{remark}

\begin{remark}
    We note that if we wish to work with a stronger vanishing Carleson measure hypothesis, (see for example Definition 1.9 of \cite{DP24}), while we have not worked out the modifications here, one should be routinely able to adopt the methods for the proof of Claim 1, to eventually pass to Carleson boxes of small enough radius for which the Carleson norm becomes arbitrarily small, and then analogously prove Claim 2 as well in this case. In this case, one doesn't need to work in `ample' subdomains as defined here.
\end{remark}

\begin{remark}
     We note that, in our construction of the `ample' sawtooth domains, we have in the end removed the Carleson boxes whose upper half consist of the 'bad' cubes that fail the condition of average smallness (ASA), and whose union has arbitrarily small measure. In the manner of the argument presented here, one cannot start with the BMO solvability assumption on the entire domain, or on any arbitrary `ample' sawtooth domains. 
\end{remark}

\section{Conclusion}

In the future, we would also like to study the converse question, of whether the natural BMO solvability estimate is implied by the weak $A_{\infty}$ condition on the elliptic measure, in presence of the Carleson measure condition on the drift. As mentioned in \cite{HL18}, this question is open even in the absence of the drift. 

\section{Notation table}
\section*{Notation table}

\subsection*{Ambient domain and operator}
\begin{longtable}{L D R}
\toprule
Symbol & Meaning & First used \\
\midrule
\endfirsthead
\toprule
Symbol & Meaning & First used \\
\midrule
\endhead
$B(0,1)\subset\mathbb{R}^{n+1}$, $n\ge 2$ & Ambient unit ball & Sec.~1 \\
$\partial\Omega = S(0,1)$ & Boundary of the unit ball & Sec.~1 \\
$L = -\operatorname{div}(A(x)\nabla(\cdot)) + B\cdot\nabla(\cdot)$ & The elliptic operator under study & eq.~(1) \\
$A(x)$ & Elliptic coefficient matrix; ellipticity constants $0<\lambda<\Lambda<\infty$ & eq.~(2) \\
$\mathcal{B}$ & Singular drift term & eq.~(1) \\
$\delta_\Omega(x) := \operatorname{dist}(x,\partial\Omega)$ & Distance from $x$ to the boundary of the domain $\Omega$ under consideration   \\
$M$ & Constant in the ambient pointwise bound $|B(x)|\le M/\delta_\Omega(x)$ & eq.~(4) \\
$\eta$ & Sawtooth parameter; fraction of $\partial\Omega$ removed & Def.~7 \\
$\Omega_\eta$ & The `ample' sawtooth domain with parameter $\eta$ (this is the object Definition~7 constructs and Theorem~2 states the result for) & Def.~7, Thm.~2 \\
\bottomrule
\end{longtable}

\subsection*{Dyadic grid and Whitney structure}
\begin{longtable}{L D R}
\toprule
Symbol & Meaning & First used \\
\midrule
\endfirsthead
\toprule
Symbol & Meaning & First used \\
\midrule
\endhead
$\mathbb{D}_k$, $Q^k_j$ & Dyadic grid on $\partial\Omega=S(0,1)$ at generation $k$; $j$-th cube in $D_k$ & Lemma~1 \\
$\sigma$ & Surface measure on $\partial\Omega$; $\sigma := \mathcal H^n|_{\partial\Omega}$ (identified with $\mathcal H^n$ on cubes) & after Lemma~1 \\
$\Delta(x,r) := B(x,r)\cap\partial\Omega$ & Surface ball & eq.~(3) \\
$\mathcal W(B(0,1))$ & Whitney decomposition of $B(0,1)$ & Sec.~1.1 \\
$U_Q$ & Whitney cube corresponding to dyadic cube $Q$ & Sec.~1.1 \\
$T_Q$ & Carleson box corresponding to $Q$; $T_Q = Q\times\{1-2^{-k}<r<1\}$ for $Q\in D_{k+1}$ & Sec.~1.1 \\
$S_Q := T_Q\setminus U_Q$ & `Bottom half' of the Carleson box $T_Q$ & eq.~(8) \\
$l(Q)$, $l(P)$, $l(U_Q)$ & Side length / diameter (up to constants) of a dyadic cube, Whitney cube, or Whitney box & Sec.~1.1 \\
$\mathcal W'(B(0,1))$ & Refined Whitney decomposition (congruent sub-cubes of $\mathcal W(B(0,1))$) & after eq.~(9) \\
$\delta_{Q,1},\dots,\delta_{Q,2(n+1)}$ & The $2(n+1)$ faces of $U_Q$; $\delta_{Q,1}$ top, $\delta_{Q,2}$ bottom, rest lateral &  \\
$\delta'_{Q,1},\dots,\delta'_{Q,2(n+1)}$ & Corresponding faces of $S_Q$ & \\
\bottomrule
\end{longtable}

\subsection*{Constants and scales}
\begin{longtable}{L D R}
\toprule
Symbol & Meaning & First used \\
\midrule
\endfirsthead
\toprule
Symbol & Meaning & First used \\
\midrule
\endhead
$\epsilon$ & Smallness parameter in the Average Smallness Assumption (ASA) & Def.~1 \\
$\overline{c}$ & Fixed constant tying $M,\epsilon$ to the scale $\tau_k$ & eq.~(6)--(7) \\
$c$ & \emph{Distinct, generic} symbol reused throughout the paper for ``some uniform constant'' (value varies occurrence to occurrence) & throughout \\
$\tau_k := \frac{\overline{c}\epsilon}{M}2^{-k}$ & Scale at which the drift enters the small-constant (Bourgain) regime & eq.~(7) \\
$l_0 = \lfloor M/\overline{c}\epsilon\rfloor$ & Integer defined via $\tau_k l_0 = 2^{-k}$ & eq.~(6)--(7) \\
$a\le 1/100$ & Uniform small aperture/geometric constant, used throughout for cone apertures, corkscrew ratios, etc. & throughout, \\
$\beta$ & Uniform lower bound constant in Claim~1, eq.~(24) & Claim 1 \\
$\iota,\alpha$ & Exponents in Claim~2 / Boundary Hölder regularity & Claim 2, Def.~9 \\
$N_0$ & Maximum number of `bad' Whitney cube families before Carleson mass is exhausted & eq.~(32) \\
\bottomrule
\end{longtable}

\subsection*{Good/bad cube bookkeeping and nested sawtooth domains (Section 2)}
\begin{longtable}{L D R}
\toprule
Symbol & Meaning & First used \\
\midrule
\endfirsthead
\toprule
Symbol & Meaning & First used \\
\midrule
\endhead
$P$ & A cube of $\mathcal W'(B(0,1))$; `good' if it satisfies ASA (Def.~1), `bad' otherwise & Def.~1 \\
$A$ & Index set of (indices $(i,k)$ of) Whitney cubes $U_{Q^k_i}$ containing at least one bad sub-cube $P\in\mathcal W'$ & Sec.~2, \\
$\mathfrak{K}_p$ & Indices in $A$ whose dyadic ancestry contains exactly $p$ elements of $A$ (stratification by ``depth'' of badness) & Sec.~2, \\
$\mathbb{F}_p := \bigsqcup_{(i,k)\in K_p} T_{Q^k_i}$ & Union of Carleson boxes at stratification level $p$ & eq.~(29) \\
$\mathbb{F}'_p := \bigsqcup_{(m,k)\in K_p} S_{Q^k_m}$ & Union of bottom-halves of Carleson boxes at level $p$ & eq.~(31) \\
$\mathbb{G}_p$ & Union of the corresponding dyadic surface cubes at level $p$ & eq.~(30) \\
$\Omega_p := B(0,1)\setminus F_p$ & Sawtooth domain after removing level-$p$ bad boxes & eq.~(34) \\
$\Lambda_p := B(0,1)\setminus F'_p$ & Sawtooth domain after removing level-$p$ bad bottom-halves & eq.~(34) \\
$\Omega_{N_0} = \Omega_\eta$ & Terminal domain of the stratification = the ample sawtooth domain & Sec.~2 \\
\bottomrule
\end{longtable}

\section{Use of LLM}
The notation table was generated using the Claude Sonnet 5 LLM model and checked by the author for accuracy. This model was also used to check for typographical errors and internal notational consistency. No LLM contributed to the mathematical content, the definitions, claims, or proofs of this paper.

\section{Data Availability Statement}
This manuscript has no associated data.

\end{document}